\documentclass{article}
\usepackage{amsmath}
 \usepackage{amsfonts}
\begin{document}
\title{{Singular Solutions of Hessian   Elliptic
Equations in Five Dimensions}} \author{{Nikolai Nadirashvili\thanks{CMI, 39, rue F. Joliot-Curie, 13453
Marseille  FRANCE, nicolas@cmi.univ-mrs.fr},\hskip .4 cm Serge
Vl\u adu\c t\thanks{IML, Luminy, case 907, 13288 Marseille Cedex
FRANCE and IITP RAS, B.Karetnyi,9, Moscow, RUSSIA,
vladut@iml.univ-mrs.fr} }}

%\date{}
\maketitle
\def\S{\mathbb{S}}
\def\Z{\mathbb{Z}}
\def\R{\mathbb{R}}
\def\C{\mathbb{C}}
\def\N{\mathbb{N}}
\def\H{\mathbb{H}}
\def\O{\mathbb{O}}
\def\tilde{\widetilde}

\def\n{\hfill\break} \def\al{\alpha} \def\be{\beta} \def\ga{\gamma} \def\Ga{\Gamma}
\def\om{\omega} \def\Om{\Omega} \def\ka{\kappa} \def\lm{\lambda} \def\Lm{\Lambda}
\def\dl{\delta} \def\Dl{\Delta} \def\vph{\varphi} \def\vep{\varepsilon} \def\th{\theta}
\def\Th{\Theta} \def\vth{\vartheta} \def\sg{\sigma} \def\Sg{\Sigma}
\def\bendproof{$\hfill \blacksquare$} \def\wendproof{$\hfill \square$}
\def\holim{\mathop{\rm holim}} \def\span{{\rm span}} \def\mod{{\rm mod}}
\def\rank{{\rm rank}} \def\bsl{{\backslash}}
\def\il{\int\limits} \def\pt{{\partial}} \def\lra{{\longrightarrow}}

\section{Introduction}
\bigskip

In this paper  we study a class of fully nonlinear second-order elliptic equations
 of the form
$$F(D^2u)=0\leqno(1.1)$$
defined in a domain of ${ \R}^n$. Here $D^2u$ denotes the
Hessian of the function $u$. We assume that
$F$ is a Lipschitz  function defined on the space $ S^2({ \R}^n)$ 
of ${n\times n}$ symmetric matrices  satisfying the uniform ellipticity condition,
 i.e. there exists a constant $C=C(F)\ge 1$ (called an {\it ellipticity
constant\/}) such that 
$$C^{-1}||N||\le F(M+N)-F(M) \le C||N||\;
\leqno(1.2)$$ 
for any non-negative definite symmetric matrix $N$; if
$F\in C^1(S^2({ \R}^n))$ then this condition is equivalent to
$$\frac{1}{ C'}|\xi|^2\le F_{u_{ij}}\xi_i\xi_j\le C' |\xi |^2\;,
\forall\xi\in { \R}^n\;.\leqno(1.2')$$
 Here, $u_{ij}$ denotes the partial derivative
$\pt^2 u/\pt x_i\pt x_j$. A function $u$ is called a {\it
classical\/} solution of (1) if $u\in C^2(\Om)$ and $u$ satisfies
(1.1).  Actually, any classical solution of (1.1) is a smooth
($C^{\alpha +3}$) solution, provided that $F$ is a smooth
$(C^\alpha )$ function of its arguments.

For a matrix $S \in   S^2({ \R}^n)$  we denote by $\lambda(S)=\{
\lambda_i : \lambda_1\leq...\leq\lambda_n\}
 \in { \R}^n$  the (ordered) set  of eigenvalues of the matrix $S$.
Equation (1.1) is called a Hessian equation  ([T1],[T2] cf. [CNS])
 if the function $F(S)$ depends only on the eigenvalues $\lambda(S)$ of the matrix $S$, i.e., if
 $$F(S)=f(\lambda(S)),$$
 for some function $f$  on ${ \R}^n$ invariant under  permutations of
 the coordinates.

 In other words the equation (1.1) is called Hessian if it is invariant under
 the action of the group
 $O(n)$ on $S^2({ \R}^n)$:
 $$\forall O\in O(n),\; F({^t O}\cdot S\cdot O)=F(S) \;.\leqno(1.3) $$
 The Hessian invariance relation (1.3) implies the following:

\medskip
 (a) $F$ is a smooth (real-analytic) function of its arguments if and only if $f$ is
a smooth (real-analytic) function.

\medskip
 (b) Inequalities (1.2) are equivalent to the inequalities
 $${\mu\over C_0} \leq { f ( \lambda_i+\mu)-f ( \lambda_i) } \leq C_0 \mu,
 \; \forall  \mu\ge 0,$$
 $\forall  i=1,...,n$, for some positive constant $C_0$.

\medskip
 (c) $F$ is a concave function if and only if $f$ is concave.

\medskip
 Well known examples of the Hessian equations are Laplace, Monge-Amp\`ere,
Bellman, Isaacs and Special Lagrangian equations.

 \medskip
 Bellman and Isaacs equations appear in the theory of controlled diffusion processes, see [F]. 
 Both are fully nonlinear uniformly elliptic equations of the form (1.1). The Bellman equation
 is concave in $D^2u \in  S^2({ \R}^n)$ variables. However, Isaacs operators are, in  general,
 neither concave nor convex. In a simple homogeneous form the Isaacs equation can be written
 as follows:
$$F(D^2u)=\sup_b \inf_a L_{ab}u =0, \leqno (1.4) $$
 where $L_{ab}$ is a family of linear uniformly elliptic operators of type 
 $$L=  \sum a_{ij} {\partial^2 \over  \partial x_i \partial x_j } \leqno (1.5)$$
  with an ellipticity  constant $C>0$ which depends on two parameters $a,b$.

\medskip
Consider the Dirichlet problem
\begin{equation}\label{dir}\begin{cases} F(D^2u, Du, u, x)=0 &\text{in}\;
 \Om \cr \quad \quad  u=\vph & \text{on}\; \pt\Om\;,\cr\end{cases}\end{equation}
 where  $\Omega \subset {\R}^n$ is a bounded domain with a smooth boundary $\partial \Omega$
and $\vph$ is a continuous function on $\pt\Om$.

We are interested in the problem of existence and regularity of
solutions to the Dirichlet problem (1.6) for Hessian equations and Isaacs equation.
The problem (1.6) has always a unique viscosity (weak)
solution for  fully nonlinear elliptic equations (not necessarily
Hessian equations). The viscosity solutions  satisfy the equation
(1.1) in a weak sense, and the best known interior regularity
([C],[CC],[T3]) for them is $C^{1+\epsilon }$ for some $\epsilon
> 0$. For more details see [CC], [CIL]. Until recently it remained
unclear whether non-smooth viscosity solutions exist. In the recent papers [NV1], [NV2], 
 [NV3], [NV4] the authors first proved   the existence   of non-classical
viscosity solutions to a fully nonlinear elliptic equation, and then of singular solutions to
Hessian  uniformly elliptic equation in all dimensions beginning from 12. Those
papers   use the functions 
$$w_{12,\delta}(x)={P_{12}(x)\over |x|^{\delta }},\;w_{24,\delta}(x)= {P_{24} (x)\over  |x|^{\delta }},\:
\delta\in [1,2[,
$$ 
with $P_{12}(x),P_{24}(x)$ being cubic forms as follows:
 $$P_{12}(x)=Re (q_1q_2q_3),\; x=(q_1,q_2,q_3)\in {\H}^3={ \R}^{12},$$
  $  {\H}$ being Hamiltonian quaternions,
$$ P_{24}(x)={Re((o_1\cdot o_2)\cdot o_3)}={Re(o_1\cdot(o_2\cdot o_3))},\; x=(o_1,o_2,o_3)\in {\O}^3={\R}^{24}
$$   $\O$ being the  algebra of Caley octonions.

Finally, 
the paper   [NTV] gives a construction  of  non-smooth viscosity solution in 5 dimensions which is order 2
 homogeneous, also for Hessian equations, the function
$$w_5(x)={P_{5} (x)\over |x|},\;  
$$
being such solution for the Cartan minimal cubic $$ P_{5}(x)=x_1^3+\frac{3 x_1}2\left(z_1^2 + z_2^2-2 z_3^2-2x_2^2\right)+\frac{3\sqrt 3}2\left(x_2z_1^2-x_2z_2^2 + 2z_1z_2z_3\right)$$ in 5 dimensions.
\bigskip

\noindent However, the methods of [NTV] does not work  for the function 
$w_{5,\delta}(x)=P_{5} (x)/ |x|^{\delta }, \;\delta>1$,   and thus does not give singular (i.e. not in $C^{1,1}$)   viscosity solutions to fully nonlinear equations in 5 dimensions.

In the present paper we fill the gap and prove

\bigskip
{\bf Theorem 1.1.}

{\it The function
$$w_{5,\delta}(x)=P_{5} (x)/ |x|^{1+\delta }, \;\delta\in [0,1[
$$ is a viscosity solution  to  a uniformly elliptic Hessian equation
  $(1.1)$ with a smooth functional $F$ in a unit ball $B\subset {\R}^{5}$  for  the isoparametric Cartan cubic form
$$ P_{5}(x)=x_1^3+\frac{3 x_1}2\left(z_1^2 + z_2^2-2 z_3^2-2x_2^2\right)+\frac{3\sqrt 3}2\left(x_2z_1^2-x_2z_2^2 + 2z_1z_2z_3\right)$$
with  $x=(x_1,x_2,z_1,z_2,z_3)$.}\bigskip

In particular  one gets  the optimality of  the interior $C^{1,\alpha}$-regularity of viscosity solutions to fully nonlinear equations   in dimensions 5 and more; note also that all previous constructions give only Lipschitz Hessian functional $F$. Let us recall that in the paper [NV5] it is proven that there is no order 2  homogenous solutions to elliptic equations in 4 dimensions which suggests  strongly that in 4 (and less) dimensions there is no homogenous non-classical solutions to uniformly elliptic equations. 
\medskip

 As in [NV3] we get  also that $w_{5,\delta}(x), \;\delta\in [0,1[$  is a viscosity solution  to a uniformly elliptic Isaacs equation:
\bigskip

{\bf Corollary 1.2.}

{\it The function
$$w_{5,\delta}(x)=P_{5} (x)/ |x|^{1+\delta }, \;\delta\in [0,1[
$$ is a viscosity solution  to  a uniformly elliptic  Isaacs equation
  $(1.4)$ in a unit ball $B\subset { \R}^{5}$.}

\medskip
 The rest of the paper is organized as follows: in Section 2 we recall some necessary preliminary results   and we prove our main results in Section 3.
The proof in Section 3 extensively uses MAPLE but is completely rigorous.

\section{Preliminary results }

\medskip
   Let $w=w_n$ be an odd homogeneous  function
of order $2-\delta,\: 0\leq \delta <1$, defined on a unit ball
$B =B_1\subset { \R}^n$ and smooth in $B \setminus\{0\}$. Then
the Hessian of $w$ is homogeneous of order $-\delta$. 

 Define the map
$$\Lambda :B  \longrightarrow \lambda (S) \in { \R}^n\; .$$
    $\lambda(S)=\{ \lambda_i :
\lambda_1\leq...\leq\lambda_n\}
 \in { \R}^n$ being  the (ordered) set  of eigenvalues of the matrix $S=D^2w$.
Denote  $\Sigma_n$  the permutation group  of    $\{ 1,...,n\}$. For
any $\sigma \in \Sigma_n$, let $T_{\sigma}$ be  the linear
transformation of ${ \R}^n$ given by $x_i \mapsto x_{\sigma(i)}, \; i=1,...,n.$

Let $a,b\in B$. Denote by $\mu_1(a,b)\leq ...\leq \mu_n(a,b)$  the eigenvalues of $(D^2w(a)-D^2w(b))$.

\bigskip\noindent {\bf  Lemma 2.1.} {\it  Assume that for a smooth function $g: U  \longrightarrow  { \R}$ where the domain $U$  contains
$$M:=\bigcup_{\sigma \in \Sigma_{n} }\   T_{\sigma }\Lambda (B)\subset { \R}^n $$ one has
$$g_{|M}=0.$$ 
 Assume also the condition $$ \min_{i=1,\ldots,5}\inf_{x\in M}\left\{\frac{\partial g}{\partial \lambda_i}(\lambda )\right\}>0.\leqno (2.1)$$ 
 Assume further  that for any $a,b\in B$ either $\mu_1(a,b)= ...= \mu_n(a,b)=0$ or 
 $$1/C\leq -\frac{\mu_1(a,b)}{\mu_n(a,b)}\leq C,\leqno (2.2)$$  where $C$ is a positive constant (may be, depending on $M,g$ but not on $a,b$). If $\delta> 0$ we assume additionally that      $w$ changes sign in $B$. Then $w$ is a viscosity solution in $B$ of a uniformly 
elliptic Hessian equation $(1.1)$ with a smooth $F$. Function $w$ is as well a solution to  a uniformly 
elliptic Isaacs equation.}\medskip

{\it Proof.} Denote for any $\theta>0$ by $K_{\theta }\subset {  \R}^n$ the cone $\{\lambda \in { \R}^n,\,
\lambda_i /|\lambda | >\theta\},\, $ and let $K^*_{\theta }$ be its dual cone. Let $x,y$ be orthogonal coordinates in ${ \R}^n$
such that $x=\lambda_1+...+\lambda_n$ and $y$ be the orthogonal complement of $x$. Denote
 by $p$ the orthogonal projection of ${ \R}^n$ on subspace $y$. Denote
 $$\Gamma =\{ g=0 \} \subset U,$$
 $$G=p(\Gamma ),$$
 $$m=p(M).$$
 From (2.1), (2.2) it follows that the surface $\Gamma $ is a graph of a smooth function $h$
 defined on $G$.  By $k_{\theta}$ we denote the function on $y$ which graph is the surface
$\partial K^*_{\theta }$.  We define the function $H(y)$ by 
$$H(y) = \inf_{z\in G}\{ h(z)+k_{\theta }(y-z)\}.$$
We fix a sufficiently small $\theta >0$. Then from (2.1), (2.2) it follows that  $H=h$ on $G$.
Denote by $J$ the graph of $H$. It is easy to show, see similar argument in [NV1], [NV3],
that for any $a,b\in J,\, a\neq b$,
$$1/C\leq -\min_i(a_i-b_i)/\max_i(a_i-b_i)\leq C.\leqno (2.3)$$
Let $E$ be a smooth function in ${ \R}^{n-1}$ with the support in a unit ball and with the integral being equal to $1$. Denote $E_c(y)=c^{-n+1}E(y/c)$, $c>0$. Set
$$H_c=H\ast E_c.$$
Then $H_c$ will be a smooth function such that any two points $a,b$ on its graph will satisfy (2.3). Moreover $H_c \rightarrow H$ in $C({ \R}^n)$ as $c$ goes to $0$, and $H_c \rightarrow h$ in $C^{\infty }$ on compact subdomains of $G$. Thus for a sufficiently small $c>0$ we can easily modify function $H_c$ to a function $\tilde H$ such that $\tilde H$ will coincide with $h$
in a neighborhood of $m$, coincide with $H$ in the complement  of $G$ and the points 
on the graph of $\tilde H$ will still satisfy (2.3) possibly with a larger constant $C$. Define the function $F$ in ${ \R}^n$ by $$F=x-\tilde H(y).$$
Then $w$ is a solution in ${ \R}^n\setminus \{ 0\}$ of a uniformly elliptic Hessian equation (1.1)
with such defined nonlinearity $F$. As in [NV3], [NV4] it follows that $w$ is a viscosity solution
of (1.1) in the whole space ${ \R}^n$. In [NV3] we have shown that the 
equation (1.1) for the function $w$ can be rewritten in the form of the Isaacs equation. The lemma is proved.\medskip

We will apply this result to the function $w_{5,\delta}(x)=P_{5} (x)/ |x|^{1+\delta } $.\medskip

Let then recall some facts from [NTV] about  the Cartan  cubic form $P_5(x)$.

\medskip

\noindent{\bf Lemma 2.2.}

 {\em The form  $P_{5}(x)$ admits a three-dimensional automorphism group.}

\medskip Indeed, one easily 
verifies that the orthogonal trasformations 

$$ A_1(\phi):= \left(%
\begin{array}{ccccc}
   {3\cos(\phi)^2-1\over 2}& { \sqrt 3\sin(\phi)^2\over 2}&0&0&{ \sqrt 3\sin(2\phi)\over 2} \\
{ \sqrt 3\sin(\phi)^2\over 2}&{1+\cos(\phi)^2\over 2}&0&\;0&{-\sin(2\phi)\over 2} \\
  0& 0&\cos(\phi)&\sin(\phi)&0 \\
   0&0&-\sin(\phi)&\cos(\phi)&0\\
  {- \sqrt 3\sin(2\phi)\over 2}&{\sin(2\phi)\over 2} &0&0&\cos(2\phi)\\
\end{array}%
\right) $$

$$ A_2(\psi):= \left(%
\begin{array}{ccccc}
   1& 0& 0& 0& 0\\
0&\cos(2\psi)& 0& -\sin(2\psi)& 0\\
0&0& \cos(\psi)&0&-\sin(\psi)\\
0&\sin(2\psi)& 0&\cos(2\psi)& 0\\
0&0&\sin(\psi)& 0& \cos(\psi)\\
\end{array}%
\right) $$

$$ A_3(\theta):= \left(%
\begin{array}{ccccc}
{3\cos(\theta)^2-1\over 2}& {-\sqrt 3\sin(\theta)^2\over 2}&0& 0&  {-\sqrt 3\sin(2\theta)\over 2} \\ 
{-\sqrt 3\sin(\theta)^2\over 2}& {1+\cos(\theta)^2\over 2}& 0& 0& {-\sin(2\theta)\over 2}\\
0&0&\cos(\theta)&-\sin(\theta)&0\\
0& 0& \sin(\theta)&\cos(\theta)& 0\\
{\sqrt 3\sin(2\theta)\over 2}&  {\sin(2\theta)\over 2}& 0& 0&\cos(2\theta)\\
\end{array}%
\right) $$
do not change the value of $P_5(x)$.
\medskip

\noindent{\bf Lemma 2.3.}

 {\em Let $G_P$ be subgroup of $SO(5)$ generated by

 \{$A_1(\phi),A_2(\psi), A_3(\theta):\:(\phi,\psi,\theta)\in { \R}^3\}.$ Then the orbit $G_PS^1$ of
 the circle
 $$S^1=\{(\cos(\chi),0,\sin(\chi),0,0):\chi\in { \R}\}\subset S^4  $$
under the natural action of  $G_P$ is the whole $S^4 . $ }
 \bigskip

We need also the following two simple algebraic results ( [NV3, Lemmas 2.2 and 4.1]):

\bigskip\noindent {\bf  Lemma 2.4.} 
  {\em Let $ A,B$  be two real symmetric  matrices
with the eigenvalues $
\lambda_1\ge\lambda_2\ge\ldots\ge\lambda_{n} $ and $
\lambda'_1\ge\lambda'_2\ge\ldots\ge\lambda'_{n} $ respectively.
Then for the eigenvalues $
\Lambda_1\ge\Lambda_2\ge\ldots\ge\Lambda_{n} $ of the matrix $A-B$
we have}
$$ \Lambda_1\ge\max_{i=1,\cdots, n}(\lambda_i-\lambda'_i),
 \;\;\Lambda_n\le\min_{i=1,\cdots, n}(\lambda_i-\lambda'_i).$$

 \bigskip\noindent
{\bf Lemma 2.5. }{\em  Let  $\delta\in[0,1),\;
     W,\;\bar W\in { \R}$ with $|W|\le \frac{1}{3\sqrt 3},|\bar W|\le \frac{1}{3\sqrt 3}$ and let
 $ \mu_1(\delta)\ge \mu_2(\delta)\ge \mu_3(\delta)$
 $($resp.,
  $\bar\mu_1(\delta)\ge\bar\mu_2(\delta)\ge\bar\mu_3(\delta)\;)$
   be the roots of the polynomial
  $$T^3+3W(1+\delta) T^2+(3W^2(1+\delta)^2-1)T+W(1-\delta) +W^3(1+\delta)^3$$
$($resp. of the polynomial
$$T^3+3\bar W (1+\delta) T^2+(3 \bar W^2(1+\delta)^2-1)T+\bar W(1-\delta)+ \bar W^3(1+\delta)^3\; ).$$ Then  for any $K>0$ verifying
$ |K-1|+|\bar W-W|\neq 0$ one has
$${1-\delta\over 5+\delta}=:\rho\le { \mu_+(K)\over -\mu_-(K)}\le
 {1\over\rho}= {5+\delta\over 1-\delta}$$
where
$$\mu_-(K):=  \min\{\mu_1(\delta)-K\bar\mu_1(\delta),\;
\mu_2(\delta)-K\bar\mu_2(\delta),\;
\mu_3(\delta)-K\bar\mu_3(\delta)\} ,$$
$$\mu_+(K):=  \max\{\mu_1(\delta)-K\bar\mu_1(\delta),\;
\mu_2(\delta)-K\bar\mu_2(\delta),\;
\mu_3(\delta)-K\bar\mu_3(\delta)\}\;.  $$   }
 
\section{Proofs}

Let $w_{5,\delta}=P_5/|x|^{1+\delta}, \: \delta\in [0,1[$. By Lemma 2.1  it is sufficient to prove the existence of a smooth function $g$ verifying  the conditions (2.1) and (2.2). We beging with calculating the eigenvalues of $D^2w_{5,\delta}(x)$.
More precisely, we need \medskip

{\bf Lemma 3.1.} 

{\em Let  $x\in S^4$, and let   $x\in G_P(p,0,q,0,0)$ with $p^2+q^2=1$.
Then $$Spec(D^2w_{5,\delta}(x))= \{\mu_{1,\delta},\mu_{2,\delta},\mu_{3,\delta},\mu_{4,\delta},\mu_{5,\delta} \}$$
for $$\mu_{1,\delta}= \frac{p(p^2\delta+6-3\delta)}{2},$$   
$$\mu_{2,\delta}= \frac{p(p^2\delta-3-3\delta)+3\sqrt{12-3p^2}}{2}, $$
$$\mu_{3,\delta}=\frac{p(p^2\delta-3-3\delta)-3\sqrt{12-3p^2}}{2}, $$ 
$$\mu_{4,\delta}= -\frac{p\delta(6-\delta)(3-p^2)+\sqrt{D(p,\delta)}}{4},$$
$$\mu_{5,\delta}= -\frac{p\delta(6-\delta)(3-p^2)-\sqrt{D(p,\delta)}}{4} ,$$
and
$$D(p,\delta):=(6-\delta)(4-\delta)(2-\delta)\delta(p^2-3)^2p^2+144(\delta-2)^2>0. $$

The characteristic polynomial $F(S)$ of $D^2w$ is given by 
$$F(S)=S^5+a_{1,\delta} S^4+a_{2,\delta}S^3+a_{3,\delta}S^2+a_{4,\delta} S+a_{5,\delta}$$
for $$a_{1,\delta}= \frac{(\delta+1)(\delta-8)b}{2},$$
$$a_{2,\delta}= \frac{(\delta+1)(21\delta+13-4\delta^2)b^2}{4}+9(2\delta-\delta^2-4),$$
$$a_{3,\delta}= \frac{(6\delta^2-31\delta-1)(\delta+1)^2b^3}{8}+
\frac{27(4\delta-2\delta^2+5+\delta^3)}{2},$$
$$a_{4,\delta}= \frac{(2\delta-1)(5-\delta)(\delta+1)^2b^4}{8}+
\frac{9(\delta-1)(\delta^2-2\delta+9)}{2},$$
$$a_{5,\delta}= \frac{b(1-\delta)\left(b^2(\delta+1)^3+ 108(1-\delta)\right)\left(b^2(\delta+1)(\delta-5)+36(\delta-1)\right)}{32},$$
$$ where \quad b:=p(p^2-3).$$}

\medskip
Note that the spectrum in this lemma is unordered one.\medskip

{\em Proof of Lemma 3.1.} Since $w_{5,\delta}$ is invariant under $G_P$, we can suppose that $x=(p,0,q,0,0)$.
Then $w_{5,\delta}(x)=\frac{p(3-p^2)}{2}$ and we get by a brute force calculation:
$$ D^2w_{5,\delta}(x):= \left(%
\begin{array}{cc} M_{1,\delta}&0\\
0&M_{2,\delta}
  \\
\end{array}%
\right) $$ being a block matrix with
$$ M_{1,\delta}:= \frac{ 1}{2}\left(%
\begin{array}{ccc}
    m_{1,1}&  {m_{1,2}}& {m_{1,3}} \\ 
{m_{1,2}}& {m_{2,2}}&  {m_{2,3}}\\
 {m_{1,3}}& {m_{2,3}}& { m_{3,3}} \\
  
\end{array}%
\right), \;$$

$m_{1,1}:=-(\delta+2)\delta p^5+(\delta+3)\delta p^3+(12-9\delta)p,$ \smallskip

 $m_{1,2}:= 3\sqrt 3 p(p^2-1)\delta,$\smallskip

 $ m_{1,3}:= -q\left((\delta+2)\delta p^4+3\delta(1-\delta)p^2+3\delta-6)\right)$ \smallskip

$ m_{2,2}:=\delta p^3-3(\delta+4)p,$ \smallskip

$ m_{2,3}:=3\sqrt 3q(\delta p^2+2-\delta),  $\smallskip

 $m_{3,3}:=(\delta+2)\delta p^5+(5-4\delta)\delta p^3-3(\delta-1)(2-\delta)p,$\smallskip
 
\bigskip

 $$
 M_{2,\delta}:=\frac{ 1}{2} \left(%
\begin{array}{cc}
{ \delta p^3+3(2-\delta)p}& 6\sqrt 3 q \\
 6 \sqrt 3 q & { \delta p^3-3(4+\delta)p}\\ 
\end{array}%
\right) $$
which gives for the characteristic polynomial $F(S) =F_1(S)\cdot F_2(S)\cdot F_3(S)$ where  
 $$F_1(S):=S-\frac{p(p^2\delta+6-3\delta)}{2}; $$
$$ F_2(S)=S^2+\frac{\delta p (p^2-3) (\delta-6) S}{2}+\frac{(2-\delta)\left( (\delta-6)\delta p^6+6(6-\delta)\delta^2 p^4+9(\delta^2-6\delta) p^2+36(\delta-2)\right)}{4};$$ 
$$F_3(S) := S^2+(3 +3\delta-\delta p^2) p S+ \frac{ (p^2-3) (\delta^2 p^4-3\delta^2 p^2-6\delta p^2+36)}{4}; $$ 
and the spectrum. Developing $F(S)$ we get the last formulas.

\medskip{\bf Corollary 3.1.} {\em Denote $\varepsilon=1- \delta$. The function $w$ verifies the following Hessian equation:
$$det(D^2w)= e_5(\Delta(w))^5  +e_3(\Delta(w))^3 S_2(w)+
 e_1\Delta(w) S_4(w)$$ where

$$e_5=\frac{ \varepsilon^2(168-5   \varepsilon^4-24   \varepsilon^3-56   \varepsilon)}{(  \varepsilon^2+3)(  \varepsilon+7)^5(  \varepsilon-2)^3} $$

$$e_3=\frac{  \varepsilon^2   (2   \varepsilon^2+  \varepsilon+8)}{ (  \varepsilon-2)^2(  \varepsilon+7)^3(  \varepsilon^2+3)},\;e_1=\frac{ \varepsilon}{ (2-  \varepsilon)( \varepsilon+7)}$$
$\Delta(w)=trace(D^2w)$ being the Laplacian, $S_2(w)$ and $S_4(w)$
being  respectively the second and the forth symmetric functions of the eigenvalues of  $D^2w$.}

\medskip
{\em Proof.} This follows immediately from Lemma 3.1 and a simple calculation since
 $$\Delta(w)=-a_{1,\delta},\;S_2(w)= a_{2,\delta},\: S_4(w)= a_{4,\delta},\:\det(D^2w)=-a_{5,\delta}.$$
\medskip

Let then  determine   the ordered spectrum $\{\lambda_1,\lambda_2,\ldots,\lambda_5 \}, \; $ $\linebreak\lambda_1\ge\lambda_2\ge\ldots\ge\lambda_5 $ of $D^2w$.  

\medskip{\bf Lemma 3.2.}  
 
{\em Let  
$\lambda_1\ge\lambda_2\ge\ldots\ge\lambda_5$ be the eigenvalues  of $D^2w_{5,\delta}(x).$
Then

$$\lambda_1= \mu_{2,\delta}, \quad \lambda_5= \mu_{3,\delta},\quad\quad\quad\quad$$
\begin{equation*}\lambda_2 = \begin{cases} \mu_{4,\delta} & \text{for}\; p\in [-1,p_0(\delta)],\cr
\mu_{1,\delta} & \text{for}\; p\in [p_0(\delta),1],  \end{cases}\end{equation*}

\begin{equation*}\lambda_3 = \begin{cases} \mu_{5,\delta} & \text{for}\; p\in [-1,-p_0(\delta)],\cr
\mu_{1,\delta} & \text{for}\; p\in [-p_0(\delta),p_0(\delta)], \cr 
\mu_{4,\delta} & \text{for}\; p\in [ p_0(\delta),1], \end{cases}\end{equation*}
\begin{equation*}\lambda_4 = \begin{cases} \mu_{1,\delta} & \text{for}\; p\in [-1,-p_0(\delta)],\cr\mu_{5,\delta} & \text{for}\; p\in [-p_0(\delta),1],\cr \end{cases}\end{equation*}

where $$p_0(\delta):=\frac{3^{1/4}\sqrt{1-\delta}}{(3+2\delta-\delta^2)^{1/4}}=\frac{3^{1/4}\sqrt{ \varepsilon}}{(4 - \varepsilon^2)^{1/4}}\in ]0,1].$$}
\medskip

{\em Proof.} The inequalities $\mu_{2,\delta}(p)\ge\mu_{1,\delta}(p)\ge\mu_{3,\delta}(p)$
are obvious since $\mu_{2,\delta}(p)$ and $\mu_{3,\delta}(p)$ are decreasing in $p$, $\mu_{1,\delta}(p)$ is increasing in $p$, $\mu_{3,\delta}(-1)=\mu_{1,\delta}(-1),$ \linebreak $\mu_{2,\delta}(1)=\mu_{1,\delta}(1).$ 

The resultant $$R(\delta,p)= Res(F_2, F_3)= 144(p-1)^2(p+1)^2\left(r_8p^8-r_6p^6+ r_4p^4-r_2p^2 +r_0  \right)$$
where
$$r_8=(\varepsilon^2-4)^2,r_6=12 (\varepsilon^2-4)^2,r_4=3 (4-\varepsilon^2)(72-17\varepsilon^2 ),$$ $$r_2=108 (\varepsilon^2-4)^2 ,r_0=144(3-\varepsilon^2)^2$$
is strictly positive for $(\varepsilon,p)\in ]0,1[\times ]-1,1[$. Indeed, let $$r:=\frac{  R}{144(p-1)^2(p+1)^2}=r_8p^8-r_6p^6+ r_4p^4-r_2p^2 +r_0$$ then
$$d:=\frac{\partial r}{4\varepsilon\partial \varepsilon}=(\varepsilon^2-4)p^8 +12p^6(4-\varepsilon^2)+3(17\varepsilon^2-70)p^4+
108(4-\varepsilon^2)p^2 +144(\varepsilon^2-3)<0$$
for $(\varepsilon,p)\in ]0,1[\times [0,1[$ since $$\frac{\partial d}{ 4p\partial p}=(4-\varepsilon^2)(-2p^6+18p^4-51p^2+54)-6p^2\ge (4-\varepsilon^2)\cdot 19-6\ge 51,$$  and for $p=1$ one has $d=-166+76\varepsilon^2\le -90.$
For $\delta=0, \varepsilon=1$ we get 
$$R(\delta,p)\ge R(1,p)=9(1-p^2)(4-p^2))(p^4-7p^2+16)$$
which proves the positivity. Using then the inequalities  
$$\mu_{2,\delta}(-1)=\mu_{4,\delta}(-1)>\mu_{5,\delta}(-1)>\mu_{3,\delta}(-1),$$
$$\mu_{2,\delta}(1)>\mu_{4,\delta}(1)>\mu_{5,\delta}(-1)=\mu_{3,\delta}(-1)$$
 and the postivity of the resultant we get
$$\mu_{2,\delta}(p)\ge\mu_{4,\delta}(p)\ge\mu_{5,\delta}(p)\ge\mu_{3,\delta}(p)$$
for any $p\in[-1,1].$

Calculatig then 
$$R_1(\delta,p)= Res(F_1, F_3)=12(p^2-3) \left( p^4(\varepsilon^2-4)+3\varepsilon^2\right)$$  
and taking into account the equalities $$\mu_{4,\delta}(p_0(\delta))=\mu_{1,\delta}(p_0(\delta)),\;\mu_{5,\delta}(-p_0(\delta))=\mu_{1,\delta}(-p_0(\delta))$$
we get the result.
\medskip

Note  the oddness property of the spectrum:
$$\lambda_{1,\delta}(-p)=-\lambda_{5,\delta}(p),\;\lambda_{2,\delta}(-p)=-\lambda_{4,\delta}(p),\;\lambda_{3,\delta}(-p)=-\lambda_{3,\delta}(p).$$

\medskip 
Let us now verify the second condition (2.2) of Lemma 2.1, namely the uniform hyperbolicity of $M_{\delta}(a,b,O)$.  
\medskip

\medskip{\bf Proposition 3.1.} 
{\em Let  $M_\delta(x)={D^2w_{\delta}}(x), \;
  0\leq\delta <1$.
 Suppose that  $a\neq b \in B_1\setminus\{0\}$  and
 let $O\in {\hbox {O}}({5} )$ be an orthogonal matrix s.t.
 $$M_\delta(a,b,O):=M_\delta(a)- {^tO}\cdot M_\delta(b)\cdot O\neq 0.$$
 Denote $ \Lambda_1\ge\Lambda_2\ge
 \ldots\ge\Lambda_{5}$  the eigenvalues of the matrix
 $M_\delta(a,b,O).$
  Then

$$   {1\over  C}\le - {\Lambda_1\over
\Lambda_{5}}\le C 
  $$
for $C:=  \frac{1000(\delta+1)(3-\delta)}{3(1-\delta)^2}.$} 
\medskip

{\em Proof.} The proof depends  on the value of 
$$ k:=p_0(\delta):=\frac{3^{1/4}\sqrt{1-\delta}}{(3+2\delta-\delta^2)^{1/4}}=\frac{3^{1/4}\sqrt{ \varepsilon}}{(4 - \varepsilon^2)^{1/4}}.$$
  Note that $C=\frac{1000(\delta+1)(3-\delta)}{3(1-\delta)^2}= \frac{1000}{k^4}$. We shall give the proof for $k\in ]0,\frac{1}{2}]$, the proof for $k\in [\frac{1}{2},1]$ is similar,  simpler and uses $C=10^4$.
 
Suppose that the conclusion does not hold, that is for some $a\neq b $  and some
  $O\in {\hbox {O}}({5} )$  one has
 $$M_\delta(a,b,O):=M_\delta(a)- {^tO}\cdot M_\delta(b)\cdot O\neq 0,$$
 but 

$$   {1\over  C}> - {\Lambda_1\over
\Lambda_{5}}\;\text{\rm or} - {\Lambda_1\over
\Lambda_{5}} >C .$$
We can suppose without loss that $|b|\le 1=|a|\in \S^4_1.$ Let $\overline{b}:=b/|b|\in \S^4_1,$ $W:=W(a),\overline W:=W(\overline b), K:= |b| ^{-1-\delta}.$ 
Note that since for any harmonic  cubic polynomial $Q(x)$ on ${\R}^n$ and any $a\in \S_1^{n-1}\subset{\R}^n$  one has
$$Tr\left (D^2\left ({Q(x)\over |x|^{1+\delta}}\right )(a)\right)= (\delta^2-2\delta-3-n) Q(a),$$
we get $Tr\left(M_{\delta}(a,b,O)\right)=(2\delta+8-\delta^2)(K\overline W-W),\; P_5   $ being harmonic. Let us prove the claim for $(K\overline W-W)\ge 0$, the proof for $(K\overline W-W)\le 0$ being the same while permuting $a$ with $b$ and $\Lambda_{1}$ with $\Lambda_{5}$.  Since  $$Tr\left(M_{\delta}(a,b,O)\right)=(2\delta+8-\delta^2)(K\overline W -W)\ge 0,$$
 we get $ 4\Lambda_{1}+\Lambda_{5}\ge 0$ and $-\Lambda_{5}/\Lambda_{1}\le 4$. Therefore, we have only to rule out the inequality $   {1\over  C}> - {\Lambda_1\over\Lambda_{5}}$ i.e. $   -\Lambda_{5}> C {\Lambda_1}.$ 
Recall that $$W=\frac{3p-p^3}{2},\overline W=\frac{3\overline p-\overline p^3}{2}$$
for some $p, \overline p\in [-1,1].$

We have then 3 possibilities:\smallskip

1). $p,\overline  p  \in [-k,k]$;\smallskip

2). $p \in [-k,k], \overline p  \notin  [-k,k]$;
\smallskip

3). $p,\overline p  \notin  [-k,k]$.\smallskip

In the cases 1) and 3) applying Lemma 2.4  we get $\Lambda_1\ge \mu_+(K)$,   $\Lambda_{5}\le \mu_-(K)$ in the notation of Lemma 2.5 which permits to finish the proof as in Proposition 4.1 of [NV3]. We thus have to treat the (more difficult)
case 2). Lemma 2.4 together with the inequality $   -\Lambda_{5}> C {\Lambda_1}$   gives $$-\min_{i=1,\cdots,5}\{K\lambda_{i,\delta}(\overline p )-\lambda_{i,\delta}( p )\}>C\max_{i=1,\cdots,5}\{K\lambda_{i,\delta}(\overline p )-\lambda_{i,\delta}( p )\}.$$

Thank to the  oddness of the spectrum we suppose without loss that 
$\overline p >k$. Recall that then by Lemma 3.2 one has  
$$ \lambda_{1,\delta}(\overline p )=\mu_{2,\delta}(\overline p ),\lambda_{1,\delta}( p )=\mu_{2,\delta}(p ),  \;\lambda_{2,\delta}(\overline p )=\mu_{1,\delta}(\overline p ),\lambda_{1,\delta}( p )=\mu_{4,\delta}(p ), $$ 
$$\lambda_{3,\delta}(\overline p )=\mu_{4,\delta}(\overline p ),\lambda_{3,\delta}( p )=\mu_{1,\delta}(p ), \;\lambda_{4,\delta}(\overline p )=\mu_{5,\delta}(\overline p ),\;\lambda_{4,\delta}( p )=\mu_{5,\delta}(p ),$$
$$\lambda_{5,\delta}(\overline p )=\mu_{3,\delta}(\overline p ),\;\lambda_{5,\delta}( p )=\mu_{3,\delta}(p ).$$

 We have then 2 possibilities for $p:$
\smallskip

2a). $p  \in [-k,0]$;\smallskip

2b). $p  \in ]0,k]$.\smallskip

Let  $p  \in [-k,0]$, then $\mu_{1,\delta}(p )\le 0$ and thus $$C\max_{i=1,\cdots,5}\{K\lambda_{i,\delta}(\overline p )-\lambda_{i,\delta}( p )\}\ge CK\lambda_{3,\delta}(\overline p ) =CK\mu_{4,\delta}(\overline p )\ge CK\mu_{4,\delta}(p_0(\delta) )=$$ $$=CK\mu_{1,\delta}(k) = CK(k^3(\sqrt {k^4+3}-k^2+3)/\sqrt {k^4+3})\ge 2CKk^{3}$$ 
since one verifies that the function $\mu_{4,\delta}(p)$ is increasing on $[k,1].$

On the other hand, $$|\min_{i=1,\cdots,5}\{K\lambda_{i,\delta}(\overline p )-\lambda_{i,\delta}( p )\}|\le K\max_{i=1,\cdots,5, p}|\{\lambda_{i,\delta}( p )\}|+\max_{i=1,\cdots,5, p}|\{\lambda_{i,\delta}( p )\}|\le 8(K+1).$$

Therefore one gets $8(K+1)\ge  2CKk^{3}$ which  clearly is   a contradiction for, say, $K\ge 1/4$. For  $0<K\le 1/4$ we get $$C\max_{i=1,\cdots,5}\{K\lambda_{i,\delta}(\overline p )-\lambda_{i,\delta}( p )\}\ge C(K\lambda_{5,\delta}(\overline p )-\lambda_{5,\delta}( p ))\ge C(K(-8)-(-5))\ge 3C$$
which can not be less than $ 8(K+1)\le 10$.

Let finally $p  \in ]0,k]$. We consider then 2 possibilities for $K$:\smallskip

$(i)$ $K \le  20/31=(1.55)^{-1} $, 

$(ii)$ $K>20/31=(1.55)^{-1}$. \smallskip

In the case $(i)$ one has $$C\max_{i=1,\cdots,5}\{K\lambda_{i,\delta}(\overline p )-\lambda_{i,\delta}( p )\}\ge C(K\lambda_{5,\delta}(\overline p )-\lambda_{5,\delta}( p ))\ge$$ $$\ge C(K\mu_{3,\delta}(1)-\mu_{3,\delta}( 0 ))\ge C(3\sqrt 3+20(\varepsilon-8)/31)>C/30>8(K+1)$$
since $\lambda_{5,\delta}( p )=\mu_{3,\delta}(p)$ is decreasing,  $\mu_{3,\delta}(0)=-3\sqrt 3,\mu_{3,\delta}(1)=\varepsilon-8.$\smallskip

We suppose then $K>20/31=(1.55)^{-1}$. Then if $p\le 3k/4$ one has $$C\max_{i=1,\cdots,5}\{K\lambda_{i,\delta}(\overline p )-\lambda_{i,\delta}( p )\}\ge CK(\lambda_{3,\delta}(\overline p )-\lambda_{3,\delta}( p )/K) \ge CK(\mu_{4,\delta}(\overline p)-\mu_{1,\delta}( p )/K)\ge$$ $$ \ge CK(\mu_{4,\delta}(k)-\mu_{1,\delta}( 3k/4 )/K)=CK(\mu_{1,\delta}(k)-\mu_{1,\delta}( 3k/4 )/K)>$$ $$> CK\mu_{1,\delta}\left( \frac{3k}{4}\right)\left(\frac{\mu_{4,\delta}(k))}{\mu_{1,\delta}\left( \frac{3k}{4}\right)}-\frac{31}{20}\right)\ge $$ $$\ge CK\frac{\mu_{1,\delta}\left( \frac{3k}{4}\right)}{ 100}  \ge CK \frac{2k^3}{100 }>20K > 8(K+1),$$
contradiction for $k\in [0,1/2]$ since $\frac{\mu_{4,\delta}(k))}{\mu_{1,\delta}\left( \frac{3k}{4}\right)}>1.56,\mu_{1,\delta}\left( \frac{3k}{4}\right)>2k^3$ there. Thus we can suppose that $p\in ]\frac{3k}{4}, k] .$ One notes then that  $\mu_{4,\delta}(p)\le \mu_{4,\delta}(k)$ for\linebreak $p\in [\frac{k}{4},1].$ This permits to rule out the case $\overline p\ge \frac{3k}{2}$. Indeed, one has in this case 
$$C\max_{i=1,\cdots,5}\{K\lambda_{i,\delta}(\overline p )-\lambda_{i,\delta}( p )\}\ge CK(\lambda_{2,\delta}(\overline p )-\lambda_{2,\delta}( p )/K) \ge CK(\mu_{1,\delta}(\overline p)-\mu_{4,\delta}( p )/K)\ge$$ $$ \ge CK(\mu_{1,\delta}(3k/2)-\mu_{1,\delta}( k)/K)=CK(\mu_{1,\delta}(3k/2)-\mu_{1,\delta}( k)/K),$$
and one gets  a contradiction as above since $$\frac{\mu_{1,\delta}\left( \frac{3k}{2}\right)}{\mu_{1,\delta}(k)}  >2 $$ for  $k\in [0,1/2]$.

The last case to rule out is thus $ K\ge 20/31, p\in [\frac{3k}{4},k], \overline p \in [k,\frac{3k}{2}]$.
Let then $$\alpha:=k-p\in [0,\frac{k}{4}]
\subset [0,\frac{1}{8}],\;\beta:=\overline p -k\in [0,\frac{k}{2}]\subset [0,\frac{1}{4}],\; a:=\max \{\alpha,\beta\}.$$

It is easy  to verify that on $\left[\frac{3k}{4},\frac{3k}{2}\right]$ one has the following inequalities:

 $$ \frac{\partial \mu_{1,\delta}(p)}{\partial p}\ge   3k^2, \;  \frac{11k^3 }{4}\ge \mu_{1,\delta}(k)\ge \frac{5k^3 }{2}\:;$$
$$ \frac{\partial \mu_{2,\delta}(p)}{\partial p}\ge -5, \;4k-5 \ge\mu_{2,\delta}(k)\ge 4k-\frac{11}{2} \:;$$
$$  \frac{\partial \mu_{3,\delta}(p)}{\partial p}\ge -\frac{ 9}{2}, \; -5-3k\ge\mu_{3,\delta}(k)\ge -\frac{11+7k  }{2}\:;$$
$$  \frac{\partial \mu_{4,\delta}(p)}{\partial p}\ge -\frac{k }{29}, \; \frac{11k^3 }{4}\ge \mu_{4,\delta}(k)\ge \frac{5k^3 }{2}\:;$$
$$ \frac{\partial \mu_{5,\delta}(p)}{\partial p}\ge 10k^2-12, \; -10k\ge\mu_{5,\delta}(k)\ge -12k \:.$$
Let then $ K\in \left[\frac{20}{31},1\right].$ Therefore,

$$C\max_{i=1,\cdots,5}\{K\lambda_{i,\delta}(\overline p )-\lambda_{i,\delta}( p )\}\ge$$ 
$$\ge C\max \{K\mu_{1,\delta}(k+\alpha )-\mu_{1,\delta}( k-\beta), K\mu_{3,\delta}(k+\alpha ) -\mu_{3,\delta}( k-\beta)\}\ge\max\left\{ M_1,M_2\right\}= $$ $$  =C\max\left\{3(K-1)k^3a+3(K+1)k^2, (1-K) (5+3k)a-\frac{ (11+7k ) (K+1)}{10}\right\} $$ 
for linear forms $M_1,M_2$ in $K$.
Note that the minimal value of $\max\{ M_1,M_2\}$ as a function of $K$ equals (recall that our $C=1000/k^4$):
$$ \frac{1500a(40-9k)}{k^2(12k^2a+18a+11k^3+20+12k)}>\frac{1250a}{k^2}>0 $$
attained for $K=K_0:= (11k^3+20+12k)/(12k^2a+18a+11k^3+20+12k)<1$.

On the other hand,
$$ -\Lambda_5\le -\min_{i=1,\cdots,5}\{K\lambda_{i,\delta}(\overline p )-\lambda_{i,\delta}( p )\}\le \max\{ l_1,l_2,l_3,l_4,l_5\}$$
for the following linear forms (in $K$)
$$l_1:=-k^2\left(\frac{11a}{4}+3k\right)K-\frac{11a}{4}k^2+3k^3,$$
$$l_2:=\left(5a+4k-\frac{11}{2}\right)K+5a+\frac{11}{2}-4k,$$
$$l_3:= (5a+5+3k)K+5a-5-3k, $$  
$$l_4:=\left(\frac{ak}{29}-\frac{11k^3}{2}\right)K+\frac{ak}{29}+\frac{11k^3}{2}, $$
$$l_5:=((12-10k^2)a+12k)K+(12-10k^2)a-12k. $$
 To refute our inequality it is sufficient to prove that $M_i(K_{j,k})> 0$ for any triple $(i,j,k)$
with $i,j\in \{1,2\},\; i\neq j, \; k\in \{1,2,3,4,5\}$ where $ l_k(K_{j,k})=M_j(K_{j,k}).$

Explicit calculations  give (for the values $m_{ijk}:=\frac{M_i(K_{j,k})}{500ak^2}$) 
$$\frac{m_{121}}{3} = \frac{(9k^4+6k^6)a+10000+5k^7+20k^4+3k^5-2250k}{k^2((3k^4+3000)a+3k^5+2750k)}>0, $$ 
$$\frac{m_{211}}{3} = \frac{(9k^4+6k^6)a+10000+5k^7+20k^4+3k^5-2250k}{(4500-3k^6)a+5000+3000k-3k^7}>0,$$ 
$$m_{122} = \frac{60000-(60k^4+90k^2)a-192k^3+66k^4-13500k-101k^2-158k^5}{k^2((6000-10k^2)a-11k^2+8k^3+5500k)}>0,$$ 
$$m_{212} = \frac{60000-(60k^4+90k^2)a-192k^3+66k^4-13500k-101k^2-158k^5}{(10k^4+9000)a+11k^4-8k^5+10000+6000k}>0,$$ 
$$m_{123} = \frac{30000-(45k^2+30k^4)a-55k^2-6750k-33k^3+30k^4-37k^5} {k^2((3000-5k^2)a+2750k-5k^2-3k^3)} >0,$$                                              $$m_{213} = \frac{30000-(45k^2+30k^4)a-55k^2-6750k-33k^3+30k^4-37k^5} {(5k^4+4500)a+5k^4+3k^5+5000+3000k} >0,$$ 
$$m_{124} = \frac{3480000-(36k^3+24k^5)a-t(k)} {k^2((348000-4k^3)a+319000k+319k^5)} >0,$$ 
  $$m_{214} = \frac{3480000-(36k^3+24k^5)a-t(k)} {(522000+4k^5)a+580000+348000k-319k^7} >0,$$    

$$m_{215} = \frac{(9k^4+30k^6-54k^2)a+15000-u(k)} {5k^4a+1500a-6k^3+1375k-6k^2a} >0,$$     
 $$m_{125} = \frac{(9k^4+30k^6-54k^2)a+15000-u(k)} {(2250+6k^4-5k^6)a+2500+1500k+6k^5 } >0,$$                                    
where  
$t(k):=783000k+ 80k^3+48k^4+44k^6+2871k^5+1914k^7<4\cdot 10^5$,\linebreak 
$u(k):=120k^2-30k^5+3375k+18k^3-
100k^4-55k^7<2000$.\medskip

Let, finally $K\ge 1$, then $$C\Lambda_1\ge C(K\mu_{1,\delta}(k+\alpha )-\mu_{1,\delta}( k-\beta))\ge \frac{5C }{2}(K-1)k^3a+3(K+1)Ck^2=$$
 $$=L_1:= \frac{500(5kK+6a-5k)}{k^2}=\frac{2500(K-1)}{k}+\frac{3000a}{k^2},$$
and
$$-\Lambda_5\le -\min_{i=2,\cdots,5}\{K\lambda_{i,\delta}(\overline p )-\lambda_{i,\delta}( p )\}\le\max\left\{L_2, L_3,L_4,L_5\right\}$$
for the following linears forms in $K$:
$$L_2:=5(K+1)a+(1-K)(5-4k) = (5a-5+4k)(K-1)+10a,$$
$$L_3:=\frac{9}{2}(K+1)a+\frac{11+7k}{2}(K-1) =\frac{9a+11+7k}{2} (K-1)+ 9a,$$
$$L_4:=(K+1)a\frac{k}{29}-\frac{5k^3}{2}(K-1)=\left(\frac{ak}{29}-\frac{5k^3}{2}\right)(K-1)+\frac{2ak}{29}, $$
$$L_5:=12k(K+1)a+\frac{9}{2}(K-1)=3\left(4ak+\frac{3}{2}\right)(K-1)+24ak.$$
One immediately sees that both the slope and the value at $K=1$ of $L_1$ are  (much)
bigger than those of $L_i, i=2,3,4,5$ which finishes the proof.\medskip

{\em Proof of Theorem 1.1.} To prove the result it is sufficient to verify the   condition (2.1)
in Lemma 2.1, namely, that the
five partial derivatives $\frac{\partial g}{\partial \lambda_i}, i=1,\ldots,5$ are strictly positive (and bounded which is automatic thank to compacity) on the
symmetrized  image $$M:=\bigcup_{\sigma \in \Sigma_{n} }\   T_{\sigma }\Lambda (B)\subset { \R}^n $$ of the unit ball under the map $\Lambda$, $$g(\lambda_1,\ldots,\lambda_5)=det(D^2w)- e_5(\Delta(w))^5  -e_3(\Delta(w))^3 S_2(w)-e_1\Delta(w) S_4(w)$$
being our equation. By homogeneity it is sufficient to show this on $M':=\Lambda (\S_1^4)$ which 
is an algebraic curve, the union of 120 curves  $T_{\sigma }\Lambda (\S_1^4)$ and that it is sufficient, by symmetry,  to verify the condition on the curve $\Lambda (\S_1^4)$ only. A brute force calculation shows then that
$$ g_1(p,\varepsilon ):=\frac{\partial g}{\partial \lambda_1}=\sum_{i=0}^{12} m_{i}p^{i}  =$$
 $$ =m_{12}p^8(p^4-12p^2+54)+m_9 b^3+m_6p^4(p^2-\frac{3}{4})+m_2+m_0, $$
\noindent with 
$$m_{12} =3(\varepsilon^4+3\varepsilon^3-20\varepsilon^2
+12\varepsilon-56)(\varepsilon-2)^2(\varepsilon+2)^2,  m_{10} =-12m_{12}, m_{8} =54m_{12},$$
$$m_{11}=m_{1}=0, m_9=3D(p,\varepsilon)(\varepsilon+7)(\varepsilon+2)(\varepsilon^2+2)(\varepsilon-2)^2,m_7=-9m_9,m_5=27m_9,$$ 
$$m_6= 108(2-\varepsilon)(3\varepsilon^7+17\varepsilon^6-54\varepsilon^5 -152\varepsilon^4+72\varepsilon^3-42\varepsilon^2+ 384\varepsilon+1344), m_3=27m_9,$$
$$m_4=-\frac{3m_6}{4},\;\;m_2=
1944 \varepsilon^2 (2-\varepsilon) (\varepsilon^2  - 7) (\varepsilon^2  +3),\;\; m_0 = -7776 \varepsilon^2  (\varepsilon^2  + 3)
 $$ for 
$D(p,\varepsilon):=\sqrt{(16-\varepsilon^2)(4-\varepsilon^2) b^2+144\varepsilon^2}, \; b=(p^2-3)p;$

$$ g_2(p,\varepsilon ):=\frac{\partial g}{\partial \lambda_2}=\sum_{i=0}^{12} n_{i}p^{i}  $$
\noindent with 
$$n_{12} =(\varepsilon+4)(\varepsilon+1)(4-\varepsilon^2)^2 ,  n_{10} =-(\varepsilon+2)(\varepsilon^4+19\varepsilon^3+86\varepsilon^2+182\varepsilon+96)(2-\varepsilon)^2,$$  $$n_{11}=n_{1}=0, \;\; n_{8} =9(\varepsilon+2)(\varepsilon^2+10\varepsilon+6)(\varepsilon^2+3\varepsilon+8)(2-\varepsilon)^2, $$ 
$$ n_9=\varepsilon(\varepsilon+7)(\varepsilon+2)(\varepsilon^2+2)(2-\varepsilon)^2 \sqrt{3(4-p^2)},\; n_7=-9n_9,\; n_5=27n_9,$$ 
 $$n_{6} =3(2-\varepsilon)(13\varepsilon^6+115\varepsilon^5+218\varepsilon^4+
170\varepsilon^3-876\varepsilon^2-2856\varepsilon-1152),$$ $$ n_4=9(\varepsilon-2)(11\varepsilon^6+62\varepsilon^5+33\varepsilon^4
-24\varepsilon^3-348\varepsilon^2-1176\varepsilon-288),$$  
$$n_3=\varepsilon^3(2-\varepsilon)(\varepsilon+7)(3\varepsilon^2+2)\sqrt{3(4-p^2)},$$
$$ n_2=108\varepsilon^2(2-\varepsilon)
(\varepsilon^2+3)(\varepsilon^2+4\varepsilon-3), n_0=1296\varepsilon^2(\varepsilon^2+3);$$

$$ g_3(p,\varepsilon ):= \frac{\partial g}{\partial \lambda_3}=\sum_{i=0}^{6} h_{2i}p^{2i}$$
with 
                                                   
  $$h_{12} = (\varepsilon + 4) (\varepsilon + 1) (4-\varepsilon^2) ^2, h_{10}  = 2 (\varepsilon + 2)(\varepsilon^4  + \varepsilon^3  - 40 \varepsilon^2 - 70\varepsilon - 48)(2-\varepsilon )^2$$
$$  h_{8} = -18 (\varepsilon + 2)(\varepsilon^4  + 4\varepsilon^3  - 19\varepsilon^2  - 28\varepsilon - 24)(2-\varepsilon )^2,
$$
 $$ h_6 = 6 (\varepsilon - 2) (7\varepsilon^6  + 37\varepsilon^5  - 136 \varepsilon^4 - 274 \varepsilon^3  + 330\varepsilon^2  + 672\varepsilon + 576),$$
$$  h_4 = 9 (\varepsilon - 2) (2 \varepsilon^6  - \varepsilon^5 + 27\varepsilon^4- 66\varepsilon^3  - 348\varepsilon^2  - 1176\varepsilon - 288),$$
$$ h_2=108\varepsilon(2-\varepsilon)(\varepsilon^3+4\varepsilon^2- 15\varepsilon-84) (\varepsilon^2+3),$$
   $$ h_0 = -1296\varepsilon(\varepsilon^2  + 3) (\varepsilon^2  + 4\varepsilon - 14);$$

$$ g_4(p,\varepsilon ):=\frac{\partial g}{\partial \lambda_4}=g_1(-p,\varepsilon ),$$
$$ g_5(p),\varepsilon :=\frac{\partial g}{\partial \lambda_5}=g_2(-p,\varepsilon ),$$
and thus we need to consider only the functions $ g_1(p,\varepsilon ), g_2(p,\varepsilon ), g_3(p,\varepsilon )$ on the set $[-1,1]\times(0,1].$ We have to prove that for any fixed $\varepsilon\in (0,1]$ they are strictly positive.

The technique of the proof is identical for all three derivatives, and we begin with
$g_3$ wich is slightly simpler since it is a polynomial in two variables.
One can rearrange it in the form $$g_3(p,\varepsilon )= g_{37}\varepsilon^7+g_{36}\varepsilon^6+g_{35}\varepsilon^5+
g_{34}\varepsilon^4+g_{33}\varepsilon^3+g_{32}\varepsilon^2+
g_{31}\varepsilon+g_{30}$$ with 
$$g_{37}=2 q^5+42 q^3-18 q^4-108 q+18 q^2,\; g_{36}=138 q^3-2 q^5-
216 q+q^6-36 q^4-45 q^2,$$
 $$g_{35}=-92 q^5+558 q^4+2160 q+5 q^6-1296-1260 q^3+261 q^2, $$ $$g_{34}=-4 q^6+5184 q-12 q^3-
5184-1080 q^2+28 q^5-36 q^4,$$
 $$g_{33}=-1944 q^2-2520 q^4-40 q^6
+14256-10692 q+5268 q^3+520 q^5,$$
 $$g_{32}=-144 q^4+72 q^3+112 q^5+17496 q-4320 q^2-16 q^6-15552,$$
 $$
g_{31}=-736 q^5+80 q^6+2304 q^4-54432 q+54432-4608 q^3
+18576 q^2,\:g_{30}=64q^2(q-3)^4\ge 0$$
for $q=p^2\in [0,1]$.
Therefore, $$g_3(p,\varepsilon )\ge \varepsilon(\bar g_{37}\varepsilon^6+\bar g_{36}\varepsilon^5+\bar g_{35}\varepsilon^4+
\bar g_{34}\varepsilon^3+ \bar g_{33}\varepsilon^2+\bar g_{32}\varepsilon+
\bar g_{31}) $$
where $\bar g_{3i}:=\min_{q\in [0,1]}g_{3i}(q)$, and an elementary calculation gives
$$\frac{g_3(p,\varepsilon )}{\varepsilon}\ge   -64\varepsilon^5-160\varepsilon^4- 5184\varepsilon^3+ 4848\varepsilon^2-15552\varepsilon+15616 >1620 $$
for $\varepsilon\in (0,\frac{9}{10}]$. For $\varepsilon\in [\frac{9}{10},1]$ we have
${g_3(p,\varepsilon )}\ge\sum_{i=0}^{6}\bar h_{2i}q^{i}$ where $\bar h_{2i}:=\min_{\varepsilon\in [\frac{9}{10},1]} h_{2i}$ and thus
$${g_3(p,\varepsilon )}\ge -736q^5+80q^6+2304q^4-54432q+54432-4608q^3+18576q^2>4840.$$
Thus, finally ${g_3(p,\varepsilon )}\ge \min \{ 1620\varepsilon,4840\}.$

The function ${g_1(p,\varepsilon )}=s_1-t_1$ with $$s_1:=(q^3-6q^2+9q) \varepsilon^7+(5q^3-30q^2+45q-72)\varepsilon^6+
(108q^2-18q^3-162q)\varepsilon^5+ $$ 
$$(-432q+288-48q^3+288q^2)
\varepsilon^4+(24q^3-144q^2+216q)\varepsilon^3+$$
$$1512\varepsilon^2+
(1152q-768q^2+128q^3)\varepsilon+448q^2(q-3)^4\ge$$
$$-72\varepsilon^6-72\varepsilon^5+96\varepsilon^4+1512\varepsilon^2\ge 1440\varepsilon^2>0,$$
and $$ t_1:=( \varepsilon^2-4)(\varepsilon+7) (\varepsilon^2+2)bD(p,\varepsilon);$$
simplifying $s_1^2-t_1^2=(s_1-t_1)(s_1+t_1)=g_1(p,\varepsilon )(s_1+t_1)$ one finds
$$ (540 q^2-1296 q-216 q^4-4 q^6+288 q^3+48 q^5) \varepsilon^{13}+$$
$$
(-1944 q^4+3024 q^3+432 q^5-7776 q+5184+2268 q^2-36 q^6) 
\varepsilon^{12}+$$
$$(-2052 q^2+3240 q^4-5328 q^3-720 q^5+60 q^6+10368 q)
\varepsilon^{11}+$$ 
$$(936 q^6+50544 q^4+55944 q^2-11232 q^5-41472-97776 q^3
+29808 q)\varepsilon^{10}$$ 
$$+(26136 q^2+120 q^6-15696 q^3+6480 q^4-24624 q
-1440 q^5) \varepsilon^9+$$
$$(57456 q^5+156816 q-492372 q^2-4788 q^6-134784
+534528 q^3-258552 q^4)\varepsilon^8+$$
$$(180576 q^2-352 q^6-313632 q+
3168 q^3+4224 q^5-19008 q^4)\varepsilon^7+$$
$$(54432 q^4-238464 q^3+
859248 q^2-1166400 q+1008 q^6+870912-12096 q^5)\varepsilon^6+$$
$$(1118592 q^3-518400 q^4-9600 q^6-1268352 q^2+736128 q +115200 q^5)
 \varepsilon^5+$$
$$(1524096 q^4+207360 q+2286144+28224 q^6-338688 q^5+
2147904 q^2-3025152 q^3) \varepsilon^4+$$
$$(10752 q^6+580608 q^4-
903168 q^3-677376 q^2+2322432 q-129024 q^5)\varepsilon^3+$$
$$(3640320 q^3-25344 q^6+304128 q^5-1368576 q^4+ 8128512 q-
7471872 q^2)\varepsilon^2+$$
$$(3096576 q^4+4644864 q^2+57344 q^6-688128 q^5-6193152 q^3) \varepsilon\ge$$
$$ 64\varepsilon^4(-10\varepsilon^9+ 18\varepsilon^8- 648\varepsilon^6- 141\varepsilon^5-2214\varepsilon^4- 2266\varepsilon^3+5760\varepsilon^2 +35721)\ge 2.2\cdot 10^6\varepsilon^4.$$
Since $s_1+t_1\le 10^6$ one gets ${g_1(p,\varepsilon )}>2\varepsilon^4.$

Similarly, ${g_2(p,\varepsilon )}=s_2-t_2$ with a polynomial $s_2\ge 3000\varepsilon^2$ and
$$t_2 =\varepsilon (\varepsilon+7) (\varepsilon-2)t(\varepsilon,q) p^3\sqrt{3(4-q)}$$
where
$$t(\varepsilon,q):=  (\varepsilon^4-2 \varepsilon^2-8) q^3+9(- \varepsilon^4+2\varepsilon^2+8) q^2
+27 (\varepsilon^4-2\varepsilon^2-8) p^2-9(3\varepsilon^2+2) \varepsilon^2 .$$
Simplifying $s_2^2-t_2^2 =g_2(p,\varepsilon )(s_2+t_2)$ one gets a polynomial $\ge$
$$(-2560\varepsilon^9-18176\varepsilon^8-325632\varepsilon^6 -1254656\varepsilon^4  +2202112\varepsilon^2+15116544) \varepsilon^4\ge 1.5\cdot 10^7 \varepsilon^4$$ and $g_2(p,\varepsilon )\ge 15\varepsilon^4$
which finishes the proof.

\bigskip
 \centerline{REFERENCES}

\medskip
 \noindent [C] L. Caffarelli,  {\it Interior a priory estimates for solutions
 of fully nonlinear equations}, Ann. Math. 130 (1989), 189--213.

\medskip
 \noindent [CC] L. Caffarelli, X. Cabre, {\it Fully Nonlinear Elliptic
Equations}, Amer. Math. Soc., Providence, R.I., 1995.

\medskip
 \noindent [CIL]  M.G. Crandall, H. Ishii, P-L. Lions, {\it User's
guide to viscosity solutions of second order partial differential
equations,} Bull. Amer. Math. Soc. (N.S.), 27(1) (1992), 1--67.

\medskip
 \noindent [CNS] L. Caffarelli, L. Nirenberg, J. Spruck, {\it The Dirichlet
 problem for nonlinear second order elliptic equations III. Functions
  of the eigenvalues of the Hessian, } Acta Math.  155 (1985),  261--301.

  \medskip
 \noindent [F] A. Friedman, {\it Differential games, } Pure and Applied Mathematics, vol. 25, John Wiley and Sons, New York, 1971.

\medskip
\noindent [NTV] N. Nadirashvili, V. Tkachev, S. Vl\u adu\c t, {\it
A non-classical Solution to Hessian Equation from Cartan Isoparametric Cubic,} Adv.   Math. 231 (2012), 1589--1597.

\medskip
\noindent [NV1] N. Nadirashvili, S. Vl\u adu\c t, {\it
Nonclassical solutions of fully nonlinear elliptic equations,}
Geom. Func. An. 17 (2007), 1283--1296.

\medskip
\noindent [NV2] N. Nadirashvili, S. Vl\u adu\c t, {\it Singular Viscosity Solutions to Fully Nonlinear Elliptic Equations},  J. Math. Pures Appl., 89 (2008), 107--113.

\medskip
\noindent [NV3] N. Nadirashvili, S. Vl\u adu\c t, {\it Octonions and Singular Solutions of Hessian  Elliptic
Equations}, Geom. Func. An. 21 (2011), 483--498.

\medskip
\noindent [NV4] N. Nadirashvili, S. Vl\u adu\c t, {\it Singular solutions of Hessian fully nonlinear elliptic equations}, Adv. Math. 228 (2011), 1718--1741.

 \medskip
\noindent [NV5] N. Nadirashvili, S. Vl\u adu\c t, {\it Homogeneous Solutions of Fully Nonlinear Elliptic Equations in Four Dimensions,} arXiv:1201.1022, to appear in Comm. Pure Appl. Math.

\medskip
\noindent [T1] N. Trudinger, {\it Weak solutions of Hessian
equations,} Comm. Part. Diff. Eq. 22 (1997),  1251--1261.

\medskip
\noindent [T2] N. Trudinger, {\it On the Dirichlet problem for
Hessian equations,} Acta Math. 175 (1995),  151--164.

\medskip
\noindent [T3] N. Trudinger, {\it H\"older gradient estimates for
fully nonlinear elliptic equations,} Proc. Roy. Soc. Edinburgh
Sect. A 108 (1988), 57--65.

 \end{document}